\documentclass[10pt]{elsarticle}

\usepackage{amssymb}
\usepackage{amsfonts}
\usepackage{verbatim}
\usepackage{xcolor}
\usepackage{amsmath}
\usepackage{amsthm}
\usepackage{extarrows}

\usepackage[a4paper,left=30mm,width=150mm,top=25mm,bottom=25mm]{geometry}

\usepackage{graphicx}
\DeclareGraphicsExtensions{.eps,.ps,.jpg,.bmp}

\newtheorem{theorem}{Theorem}
\newtheorem*{thm}{Main Theorem}
\newtheorem{lemma}{Lemma}

\newtheorem*{remark}{Remark}
\usepackage{lineno,hyperref}
\modulolinenumbers[5]

\journal{JMP}
\begin{document}

\begin{frontmatter}

\title{Mixing for acoustic wave motion with boundary random force\tnoteref{mytitlenote}}
\tnotetext[mytitlenote]{This research was partially supported by the National Natural Science Foundation of China (12272297) and the Fundamental Research Funds for the Central Universities.}

\author[mymainaddress]{Zhe Jiao\corref{mycorrespondingauthor}}
\cortext[mycorrespondingauthor]{Corresponding author}
\ead{zjiao@nwpu.edu.cn}

\author[mymainaddress]{Xiao Li}

%\address[mymainaddress]{MOE Key Laboratory for Complexity Science in Aerospace, School of Mathematics and Statistics, Northwestern Polytechnical University, Xi'an 710129, China}

\address[mymainaddress]{School of Mathematics and Statistics, Northwestern Polytechnical University, Xi'an 710129, China}

\begin{abstract}
This paper concerns about the large time behavior of acoustic wave motion driven by a random force acting through the boundary. 
We begin with an abstract result showing the interconnection between the regularity of Markov semigroup generated by a stochastic evolution equation and the observability property of the corresponding adjoint system. This result is then applied to study the mixing for acoustic wave system with a boundary random perturbation of the white noise type.  
We shall show there exists a unique invariant measure for the stochastic wave system, and the law of the solution to the system converges to this invariant measure weakly.  
\end{abstract}


\end{frontmatter}

%\linenumbers

\section{Introduction}
\label{sec:intro}
Let $D\subset \mathbb{R}^3$ be an open and bounded domain that is full of some kind of idealized fluid. The equations of acoustic wave motion in $D$ will be considered below
\begin{equation}
	\rho\frac{\partial \mathbf{u}}{\partial t} = -\nabla p, \quad \kappa\frac{\partial p}{\partial t} = -\nabla\cdot\mathbf{u}, \label{A1}
\end{equation}
where the vector $\mathbf{u}(t, \mathbf{x})$ is the fluid velocity, $p(t, \mathbf{x})$ the acoustic pressure, $\rho$ the uniform density of the fluid and $\kappa$ the adiabatic compressibility, see in \cite{MI}. When the acoustic wave is incident on the smooth boundary of $D$, denoted by $\partial D$, the acoustic pressure $p$ tends to make the boundary move, and forces some fluid into the pores of the boundary. 
As for any fluid motion that is normal to the boundary, there will be wave motion in the material forming the boundary. Then on the boundary we have the following physical relationship which is called the acoustic boundary condition.
\begin{equation} 
\left
    \{
        \begin{array}{ll}
            -p=\textrm{Vibrating motion of the boundary,}\\ 
             \mathbf{n}\cdot \mathbf{u}=\textrm{Velocity of the vibrating boundary.} \label{A2}
        \end{array} 
\right.
\end{equation}
Here, $\mathbf{n}$ is the unit exterior normal vector of the boundary. And if the motion of the boundary is perturbed by a small random fluctuation, the acoustic wave motion in $D$ will be also subject to the effect of the small random perturbation. It is interesting and important to consider the random fields governed by the stochastic acoustic wave equation driven by a small random perturbation.

Since $\mathbf{u}$ can be obtained by taking the gradient of a scalar function $\phi(t, \mathbf{x})$, $\mathbf{u}=-\nabla\phi$, then $\phi$ is always called the velocity potential. And it implies from (\ref{A1}) that $p=\rho\frac{\partial\phi}{\partial t}$ and $\phi$ satisfies the following wave equation
\begin{equation}
	\frac{\partial^2 \phi}{\partial t^2} = c^2\triangle \phi, \quad \textrm{in $\mathbb{R}^{+} \times D$.} \label{1}
\end{equation}
Here, $c=\sqrt{\frac{1}{\kappa\rho}}$ is the speed of wave propagation. 
Assume that a portion $\Gamma_0$ of boundary $\partial D$ has a very low acoustic impedance compared to the acoustic impedance of the idealized flow. In this case the portion of boundary is called sound soft \cite{Nail}. Then we have
\begin{equation} %\tag{B1}
	 \phi = 0,  \quad \textrm{on $\mathbb{R}^{+} \times \Gamma_0$}. \label{2} 
\end{equation}
As for the other portion $\Gamma_{1}=\partial D-\Gamma_0$, the motion of the portion $\Gamma_1$ at one point will be related to motion at another point on $\Gamma_1$, and the relationship depends on the motion inside the material. Let $\delta(t, \mathbf{x})$ be the normal displacement into the domain of a point $\mathbf{x }\in \Gamma_1$ at time $t$. Suppose that the material forming the boundary is elastic. If the motion on the different parts of $\Gamma_1$ does not influence each other directly, each point on $\Gamma_1$ acts like a resistive harmonic oscillator, which can deduces from (\ref{A2}) the locally reacting boundary condition as follows. 
\begin{equation} %\tag{B2-1}
\left
    \{
        \begin{array}{ll}
            -\rho\phi_{t}= m\delta_{tt} +d\delta_t+k\delta+\xi(t),\\ 
            \frac{\partial\phi}{\partial \mathbf{n}}= \delta_{t},  \label{3}
        \end{array} 
\right. \quad \textrm{on $\mathbb{R}^{+} \times \Gamma_1$}.
\end{equation}
Here, $m$ is the mass per unit area of $\Gamma_1$, $d$ is the resistivity, $k$ is the stiffness, and $\xi$ is a random perturbation of the white noise type. 
The acoustic wave system (\ref{1}-\ref{2}-\ref{3}) is supplemented with the initial condtion $\phi(0)=\phi_0$, $\phi_t(0)=\phi_1$, $\delta(0)=\delta_0$ and $\delta_t(0)=\delta_1$.
 
When the boundary random force equals to zero, that is $\xi(t)=0$, this determined form of (\ref{3}) has been obtained in~\cite{BR}. The mechanical energy associated with the acoustic wave system (\ref{1}-\ref{2}-\ref{3}) is defined as $\mathcal{E}(t) = \mathcal{E}_1(t) + \mathcal{E}_2(t)$ with
\begin{equation*}
\begin{aligned}
	\mathcal{E}_1(t) =\frac{1}{2}\int_{D}\rho[c^{-2}|\phi_t|^{2} +|\nabla\phi|^{2}]dx, \quad \mathcal{E}_2(t) = \frac{1}{2}\int_{\Gamma_{1}}[k|\delta|^{2}+m|\delta_t|^{2}]d\Gamma.
\end{aligned}
\end{equation*}
The first part $\mathcal{E}_1(t)$ of the energy is generated by the motion of acoustic wave, and the second part $\mathcal{E}_2(t)$ generated by the boundary vibration.
Thus, it is not hard to check that
\begin{equation} \label{dissipation}
	\frac{d}{dt}\mathcal{E}(t) = -d\int_{\Gamma_1}|\delta_t|^2d\Gamma \leq 0,
\end{equation} 
from which we can see that the energy dissipated by the boundary friction. The author in \cite{B1} point out the energy $\mathcal{E}(t)$ have no uniform decay rates. Afterwards, authors in \cite{Qin} obtained polynomial decay of $\mathcal{E}(t)$ under the following geometrical conditions
\begin{equation} \label{G1}
\begin{split}
	\Gamma_0 =\{x\in\Gamma: (x-x_0)\cdot\mathbf{n}\leq 0 \},\quad \Gamma_1 =\{x\in\Gamma: (x-x_0)\cdot\mathbf{n}\geq c \}
\end{split}
\end{equation}
for some point $x_0\in\mathbb{R}^3$ and some constant $c>0$. We refer the reader to \cite{LL, AN1, AN2, Gao2} for deep studies on stability for the related acoustic wave systems.

However, if the boundary random force $\xi\neq 0$, what is the propagation of acoustic wave? As we know, the friction tends to drive any system to a completely \textit{dead} state, while the noise tends to keep the system \textit{alive}. What is the balance between friction and noise? From the viewpoint of mathematics, we should concern about the large time behavior of solution to the acoustic wave system~(\ref{1}-\ref{2}-\ref{3}) under the  interaction of the friction and random force both acting on the boundary. 
To be precise, we are interested in the problem of mixing for the corresponding random flow generated by the acoustic wave system~(\ref{1}-\ref{2}-\ref{3}).
The following theorem is a simplified version of the main result of this paper. (see Section~\ref{sec:mixing} for an exact statement).
\begin{thm}
Under the above geometric condition of the domain, the acoustic wave system~(\ref{1}-\ref{2}-\ref{3}) is mixing which means the law of its solution converges to an invariant measure weakly as time goes to infinity. 
\end{thm}

Let us mention that there have been relatively few mathematical analysis of the wave equation with stochastic boundary values. In 1993, Mao and Markus in \cite{Mao} used the parallelogram identity to represent the solution explicitly. However, the one-dimensionality of space is essential in their analysis. Afterwards, Kim in \cite{Kim} combined the Galerkin method and the duality argument to establish the well-posedness for a one-dimensional wave equation with variable coefficients and white noise on the Neumann boundary. On the other hand,  the parabolic equations with boundary white noises have been studied by some authors (e.g. \cite{Sowers, Maslowski, FW}). Recently, Shirikyan in \cite{Armen2} developed a general framework for dealing with random perturbations acting through the boundary of the domain and applied this result to study the two-dimensional Navier-Stokes system driven by a random force on the Dirichlet boundary. To the best of my knowledge, the problem of mixing for the acoustic wave system with a random perturbation acting through the boundary has not been studied in earlier works.

In the present paper, the acoustic wave system (\ref{1}-\ref{2}-\ref{3}) is reformulated to be a stochastic evolution equation on a Hilbert space.  We prove some properties of the corresponding operator, which implies the existence and uniqueness of invariant measures for this equation.  
By the controllability of an associated deterministic system, we can transform the mixing problem for stochastic wave system into a problem of observability. Therefore, the core of our work is to prove the observability inequality of an adjoint system.
Note that the geometry condition of the domain plays a decisive role in the proof of the observability inequality.
Our result is new with respect to the literature on two accounts: (i) stochasticity is first accounted for in the acoustic wave model; (ii) the observability of an adjoint system implies the mixing. Finally, we would like to point out that our idea is stimulated by these significant papers \cite{Kim, Gao1, Armen1, Armen2}.

This paper is organized as follows. In section \ref{sec:regularity}, we study an abstract stochastic evolution equation and prove a result on the equivalence relation between regularity of the Markov semigroup and observability of a corresponding adjoint system, which is a key lemma used in this work. The main result of this paper on the mixing for the acoustic wave system perturbed by a random boundary force is proposed in Section \ref{sec:mixing}. Section \ref{sec:conclusion} provides some conclusions. 

%For the related work, we refer to \cite{MI, JX, JYX, DFMP1, LL, JXZ}.
\subsection*{Notation}  

Let $\mathcal{H}$ and $\mathcal{L}$ be separable Hilbert spaces. We shall use the following notations. $\mathcal{B}(\mathcal{H})$ denotes the Borel $\sigma$-algebra over $\mathcal{H}$. $\mathcal{P}(\mathcal{H})$ is the set of probability Borel measure on $\mathcal{H}$. $B_b(\mathcal{H})$ ($C(\mathcal{H})$, $C_b(\mathcal{H})$, respectively) is the set of bounded Borel functions (continuous, bounded and continuous functions, respectively) on $\mathcal{H}$. $B(\mathcal{L})$ is the set of all linear bounded operators on $\mathcal{L}$.

A symmetric nonnegative operator $\mathrm{Q}\in B(\mathcal{L})$ is of trace class. There exits a completely, orthonormal basis $\{e_k\}$ in $\mathcal{L}$ and a bounded sequence $\{\gamma_k\}$ of nonnegative real number such that
\[
	\mathrm{Q}e_k = \gamma_k e_k, \quad k=1, 2, \cdots.
\]

Throughout the paper, $(\Omega, \mathcal{F}, \mathbb{P})$ is a probability space with a right-continuous increasing family $\mathcal{F}=(\mathcal{F}_t)_{t\geqslant 0}$ of sub-$\sigma$-fields of $\mathcal{F}_0$ each containing $\mathbb{P}$-null sets. $\mathbb{E}(\cdot)$ stands for expectation with respect to the probability measure $\mathbb{P}$.  

Denote by $\mathrm{Dom}(\cdot)$ the domain of an operator, and $\mathrm{Ran}(\cdot)$ the range of an operator. $C(\cdot)$ denote positive numbers, which may depend on the quantities mentioned in the brackets.

\section{Regularity in terms of observability inequality}
\label{sec:regularity}

Let $W(t)$, $t\geq 0$ be a $\mathcal{L}$-valued Q-Wiener process, and has the following expansion
\begin{equation*}
	W(t) = \sum_{k=1}^{\infty}\gamma_k \beta_{k}(t)e_{k}, \quad t\geq 0,
\end{equation*} 
where $\{\beta_{k}\}$ is a sequence of real valued Brownian motions mutually independent on the above probability space. 
In this section, we are concerned with the linear stochastic evolution equation driven by this Wiener process 
\begin{equation} \label{linearSDE}
\left
    \{
        \begin{array}{ll}
            d\mathbf{X}(t)=A(\mathbf{X}(t))dt+ BdW(t),\\
            \mathbf{X}(0)=\mathbf{x},    
        \end{array}
\right.
\end{equation}
under the following assumptions.
\begin{itemize}
\item $A$ is the infinitesimal generator of a strongly continuous semigroup $S(t)$, $t\geqslant 0$, on $\mathcal{H}$.
\item $B$ is a linear continuous mapping from $\mathcal{L}$ into $\mathcal{H}$.
\item The linear operator $Q_t$, $t > 0$, defined by
\[
	Q_{t}\mathbf{x}=\int_{0}^{t}S(s)BB^{\ast}S^{\ast}(s)\mathbf{x} ds, \quad \mathbf{x}\in\mathcal{H}
\]
is of trace class.
\end{itemize}
Here, we denote by $S^{\ast}(t)$ the adjoint of $S(t)$ and by $B^{\ast}$ the adjoint of $B$. We know $S^{\ast}(t)$ is a strongly continuous semigroup of continuous linear operators and the infinitesimal generator of this semigroup is the adjoint $A^\ast$ of $A$ (see \cite[Corollary 10.6]{Pazy}).

By Theorem 5.3.1 in \cite{PratoZabczyk}, there exists an $\mathcal{H}$-valued $\mathcal{F}_{t}$-adapted process $\mathbf{X}(t)$, $t\geq 0$, satisfies the following integral equation  
\begin{eqnarray*}
	\mathbf{X}(t) =S(t)\mathbf{x} +\int_{0}^{t}S(t-s)BdW(s), \quad t\in [0, T]
\end{eqnarray*}
for any $\mathbf{x} \in \mathcal{H}$ is $\mathcal{F}_0$-measurable and $\mathbb{E}\|\mathbf{x}\|_{\mathcal{H}}^2\leq \infty$, which means $\mathbf{X}(t)$ is a mild solution of \eqref{linearSDE}.

Let $P_t$ stands for the transition function of the family defined as the law of the solution $X(t, \mathbf{x})$ under the probability measure $\mathbb{P}$
\[
	P_t(\mathbf{x}, \Lambda)=\mathbb{P}(X(t, \mathbf{x})\in \Lambda), \quad \Lambda\in \mathcal{B}(\mathcal{H}), \quad t\geqslant0.	
\]
The corresponding Markov semigroups are given by
\begin{equation*}
 \begin{aligned}
 	M_{t}&:  C(\mathcal{H})\rightarrow C(\mathcal{H}), \quad M_{t}f(\mathbf{x})=\int_{\mathcal{H}}P_{t}(\mathbf{x}, dz)f(z),\\
	M^{\ast}_{t}&:  \mathcal{P}(\mathcal{H})\rightarrow \mathcal{P}(\mathcal{H}), \quad M^{\ast}_{t}\lambda(\Lambda)=\int_{\mathcal{H}}P_{t}(\mathbf{x}, \Lambda)\lambda(d\mathbf{x}).
\end{aligned}
\end{equation*}
Note that these two semigroups satisfy the duality relation $(M_{t}f, \lambda)=(f, M_{t}^{\ast}\lambda)$ for any $f\in C(\mathcal{H})$, $\lambda\in \mathcal{P}(\mathcal{H})$. 
A probability measure $\mu\in\mathcal{P}(\mathcal{H})$ is said to be invariant with respect to $M_{t}$ if and only if $M^{\ast}_{t}\mu = \mu$ for each $t\geqslant0$. A Markov semigroup $M_{t}$ is $t_0$-regular if all transition probabilities $P_{t_0}(\mathbf{x}, \cdot)$, $\mathbf{x}\in\mathcal{H}$, are mutually equivalent. 

In the following, we give a sufficient and necessary condition for the regularity of the Markov semigroup $M_{t}$.

\begin{theorem}\label{regularity}
$M_{t}$, $t\geqslant0$, is $t$-regular if and only if there exists a positive constant $C_{\textrm{obs}}>0$ such that
\begin{eqnarray}\label{obsInquality}
	\int_{0}^{T}\|B^{\ast}Z(t)\|^2_{\mathcal{L}}ds \geq C_{\textrm{obs}}\|Z(0)\|^2_{\mathcal{H}},\quad \forall z\in\mathrm{Dom}(A^{\ast}),
\end{eqnarray}
in which $Z(t)$ is the solution of the adjoint system
\begin{equation*} 
\left
    \{
        \begin{array}{ll}
            \dot{Z}(t)= -A^{\ast}(Z(t)),\quad t\in[0, T],\\
            Z(T) = z\in \mathrm{Dom}(A^{\ast}).   
        \end{array} \label{adjoint}
\right.
\end{equation*} 
\end{theorem}
\begin{proof}
Note that $X(t, \mathbf{x})$ is a Gaussian random variable with mean $S(t)\mathbf{x}$ and covariance $Q_{t}$. If the range of $S(t)$ is a subset of the range of $Q_{t}^{\frac{1}{2}}$, that is,  
$\textrm{Ran} S(t)\subset \textrm{Ran}Q_{t}^{\frac{1}{2}}$, $t>0$.
Then we can define a linear bounded operator $R(t)=Q_{t}^{-\frac{1}{2}}\circ S(t)$.  By the Cameron-Martin formula, we have for any $\mathbf{x}$, $\mathbf{y}\in \mathcal{H}$
 \begin{equation}\label{CMformula}
 \begin{aligned}
	P_t(\mathbf{x}, \Lambda) &= \int_{\Lambda}F(t,\mathbf{x}, \tilde{\mathbf{x}})P_t(0, d\tilde{\mathbf{x}})\\
	&= \int_{\Lambda}\frac{F(t,\mathbf{x}, \tilde{\mathbf{x}})}{F(t,\mathbf{y}, \tilde{\mathbf{x}})}F(t,\mathbf{y}, \tilde{\mathbf{x}})P_t(0, d\tilde{\mathbf{x}}) = \int_{\Lambda}\frac{F(t,\mathbf{x}, \tilde{\mathbf{x}})}{F(t,\mathbf{y}, \tilde{\mathbf{x}})}P_t(\mathbf{y}, d\tilde{\mathbf{x}})
 \end{aligned}
\end{equation} 
with
$F(t, \cdot, \tilde{\mathbf{x}})=\exp\{\big(R(t)\cdot, Q_{t}^{-\frac{1}{2}}\tilde{\mathbf{x}}\big)_{\mathcal{H}}-\frac{1}{2}\|R(t)\cdot\|_{\mathcal{H}}^2\}$.	
From \eqref{CMformula} it implies that all transition probabilities $P_{t}(\mathbf{x}, \cdot)$, $\mathbf{x}\in\mathcal{H}$, are mutually equivalent. Therefore, we have the sufficient and necessary condition for the regularity of $M_{t}$ is that $\textrm{Ran} S(t)\subset \textrm{Ran}Q_{t}^{\frac{1}{2}}$.

Next, we shall show that the condition $\textrm{Ran} S(t)\subset \textrm{Ran}Q_{t}^{\frac{1}{2}}$ is equivalent to the fact that the following controlled system
\begin{equation} 
\left
    \{
        \begin{array}{ll}
            \dot{\mathbf{Y}}(t)=A(\mathbf{Y}(t))+ B\mathbf{v}, \quad t\in[0, T],\\
            \mathbf{Y}(0)=\mathbf{y}_0 \in \mathcal{H},  
        \end{array} \label{deterministic}
\right.
\end{equation}
where $T>0$ is fixed and $\mathbf{v}\in L^2([0. T]; \mathcal{L})$, is null controllable in time $T$, which means given any $\mathbf{y}_0\in \mathcal{H}$, there exists a control function $\mathbf{v}$ such that $\mathbf{Y} (T)=0$. 

The solution of (\ref{deterministic}) can be written as
\begin{eqnarray*}
	\mathbf{Y}(t) =S(t)\mathbf{y}_0 +\int_{0}^{t}S(t-s)B\mathbf{v}ds, \quad t\in [0, T].
\end{eqnarray*}
Define the following operator
\[
	F_{T}: L^{2}([0, T]; \mathcal{L})\rightarrow \mathcal{H}, \quad F_{T}\mathbf{v}=\int_{0}^{T}S(T-s)B\mathbf{v}ds.
\]
Note that Ran$F_T$ consists of all states reachable in time $T$ from zero. Let $z\in\textrm{Dom}(A^{\ast})$. Let us also recall that $\textrm{Dom}(A^{\ast})$ is dense in $\mathcal{H}$. Then for $\mathbf{v}\in L^{2}([0, T]; \mathcal{L})$ we have
\begin{eqnarray*}
\begin{aligned}
	(\mathbf{v}, F^{\ast}_{T}(z))_{L^{2}([0, T]; \mathcal{L})}=(F_{T}(\mathbf{v}), z)_{\mathcal{H}}=\int_{0}^{T}(\mathbf{v}, B^{\ast}S^{\ast}(T-s)z)_{\mathcal{H}}ds,
\end{aligned}
\end{eqnarray*}
which implies that
\[
	(F^{\ast}_{T}\mathbf{x})(s)=B^{\ast}S^{\ast}(T-s)\mathbf{x}, \quad \mathbf{x}\in\mathcal{H}, s\in[0, T].
\]
By the definition of $Q_{T}^{\frac{1}{2}}$, we see that $Q_{T}=F_{F}F^{\ast}_{T}$, and $\textrm{Ran}Q_{T}^{\frac{1}{2}} =\textrm{Ran}F_T$.
Moreover, thanks to the representation $\mathbf{Y}(T) =S(t)\mathbf{y}_0 +F_T\mathbf{v}$, then given any $\mathbf{y}_0\in \mathcal{H}$, there exists a control function $\mathbf{v}\in L^2([0. T]; \mathcal{L})$ such that $\mathbf{Y} (T)=0$, if and only if $\textrm{Ran}S(t) \subset  \textrm{Ran}F_T = \textrm{Ran}Q_{T}^{\frac{1}{2}} $.
%Therefore, to prove $\textrm{Ran} S(t)\subset \textrm{Ran}Q_{t}^{\frac{1}{2}}$, it remains to show that system (\ref{deterministic}) is null controllable in time $T$.

From \cite[Theorem 2.44]{Coron}, the system (\ref{deterministic}) is null controllable in time $T$ if and only if there exists a positive constant $C_{\textrm{obs}}>0$ such that
\begin{equation}\label{obs1}
	\int_{0}^{T}\|B^{\ast}S^{\ast}(s)z\|^2_{\mathcal{L}}ds \geq C_{\textrm{obs}}\|S^{\ast}(T)z\|^2_{\mathcal{H}},\quad \forall z\in\textrm{Dom}(A^{\ast}).
\end{equation}
Due to $Z(t) = S^{\ast}(T-t)z$, $t\in[0, T]$, we have $Z(0) = S^{\ast}(T)z$ and
\begin{equation*}
	\int_{0}^{T}\|B^{\ast}S^{\ast}(s)z\|^2_{\mathcal{L}}ds \xlongequal{s = T-t} \int_{0}^{T}\|B^{\ast}S^{\ast}(T-t)z\|^2_{\mathcal{L}}dt = \int_{0}^{T}\|B^{\ast}Z(t)\|^2_{\mathcal{L}}dt,
\end{equation*}
which implies the inequality \eqref{obs1} is equivalent to \eqref{obsInquality}. 
\end{proof}

\begin{remark}
As usual, the inequality \eqref{obsInquality} is called \emph{observability inequality}. Due to this theorem, we can study the mixing property of a stochastic system via the proof of observability inequality for the corresponding adjoint system, which we shall see in the following section.
\end{remark}

\section{Mixing for random acoustic wave}
\label{sec:mixing}
In this section, we apply Theorem \ref{regularity} to the acoustic wave system (\ref{1}-\ref{2}-\ref{3}). We first formulate this system into the form of the stochastic evolution equation \eqref{linearSDE}.
Let $H_{0}(D)=\{v(x)\in H^1(D): v(x)=0, \textrm{on $\Gamma_0$}\}$. Define the finite energy space by
\[
	\mathcal{H}=H_{0}(D)\times L^2(D)\times L^2(\Gamma_1)\times L^2(\Gamma_1),
\]
where the inner product in $\mathcal{H}$ is given by
\[
	(\mathbf{ f}, \mathbf{ g})_{\mathcal{H}}=\int_{D}\rho[c^{-2}f_2\overline{g_{2}} +\nabla f_1\overline{\nabla g_1}]dx+ \int_{\Gamma_{1}}[kf_3\overline{g_{3}}+mf_4\overline{g_{4}}]d\Gamma.
\]
%The condition $\phi|_{\Gamma_0}=0$ and $\frac{\partial \phi}{\partial x}|_{\Gamma_1}=\delta_t$ are interpreted in the weak sense
%\begin{eqnarray*}
%	\int_{D}[\Delta \phi \overline{v} +\nabla\phi\cdot \overline{\nabla v}]dxdy =\int_{\Gamma_1}\delta_t \overline{v} d\Gamma,  \quad \forall v\in H^1(D),
%\end{eqnarray*}
%which can be deduced from integrating by parts.
Let $\mathcal{L}=L^2(D)\times L^2(D)\times L^2(\Gamma_1)\times L^2(\Gamma_1)$. Then we have the Gelfand triple
\[
	\mathcal{H}\hookrightarrow\mathcal{L}\equiv\mathcal{L}^{\ast}\hookrightarrow\mathcal{H}^\ast.
\]
Set $\mathbf{X}(t)=\{\phi, \phi_t, \delta, \delta_t\}$ and $\mathbf{x}=\{\phi_0, \phi_1, \delta_0, \delta_1\}$. The random perturbation $\xi(t)$ of the white noise type is the formal time derivative of $W(t)$. 
Then the acoustic wave system (\ref{1}-\ref{2}-\ref{3}) can be rewritten as a stochastic evolution equation \eqref{linearSDE} with
\begin{equation*}       
A=
\left(                 
  \begin{array}{cccc}   
    0 & 1 & 0 & 0\\  
    c^2\Delta & 0 & 0 & 0\\  
     0 & 0 & 0 & 1\\
     0 & -\frac{\rho}{m} & -\frac{k}{m}  & -\frac{d}{m} \\
  \end{array}
\right),                
\quad    
B=
\left(                 
  \begin{array}{cccc}   
    0 & 0 & 0 & 0\\  
    0 & 0 & 0 & 0\\  
    0 & 0 & 0 & 0\\
    0 & 0 & 0  & -\frac{1}{m} \\
  \end{array}\right)              
\end{equation*}
and the domain of the operator $A$ given by 
\[
	\mathrm{Dom}(A)=\{\mathbf{f}\in \mathcal{H}: \Delta f_1\in L^2, f_2 \in H_{\Gamma_0}, \frac{\partial f_1}{\partial n}|_{\Gamma_1}=f_4\}.
\]

%The theorem below ensure the existence and uniqueness of the solutions of the system (\ref{1}-\ref{2}-\ref{3}).
\begin{lemma}\label{operator}
(i) $A$ is a densely defined closed linear operator, which is the infinitesimal generator of a $C_0$ semigroup of contractions $S(t)$, $t\geqslant 0$.

(ii) The corresponding linear operator $Q_t$ is of trace class, and $\sup\limits_{t\geq0}\mathrm{Tr}Q_t < \infty$.

(iii) If the initial data $\mathbf{x}$ is $\mathcal{F}_0$-measurable and $\mathbb{E}\|\mathbf{x}\|_{\mathcal{H}}^2\leq \infty$, then system (\ref{1}-\ref{2}-\ref{3}) has a unique mild solution. 
\end{lemma}
\begin{proof}
(i) Due to $\mathrm{Dom}(A)$ is dense in $\mathcal{H}$, then the operator $A$ is a densely defined. Then we can deduce uniquely the duality operator $A^{\ast}$ of $A$
\begin{equation*}  
%A^{\ast}Z(t)=
%\left(                 
%  \begin{array}{c}   
%    -Z_{2} \\  
%    -c^2\Delta Z_{1}\\  
%     -Z_{4}\\
%     \frac{\rho}{m}Z_{2} +\frac{k}{m}Z_{3}  -\frac{d}{m}Z_{4} \\
%  \end{array}
%\right), 
%\quad     
A^{\ast}=
\left(                 
  \begin{array}{cccc}   
    0 & -1 & 0 & 0\\  
    -c^2\Delta & 0 & 0 & 0\\  
     0 & 0 & 0 & -1\\
     0 & \frac{\rho}{m} & \frac{k}{m}  & -\frac{d}{m} \\
  \end{array}
\right)           
\end{equation*}
and the domain of the operator $A^{\ast}$ is given by 
\[
	\mathrm{Dom}(A^{\ast})=\{\mathbf{g}\in \mathcal{H}: \Delta g_1\in L^2, g_2 \in H_{\Gamma_0}, \frac{\partial g_1}{\partial n}|_{\Gamma_1}=g_4\}.
\]
Now, we shall show that $A$ is closable. Suppose it is not true. Then there is a sequence $\mathbf{f}_n\in\mathrm{Dom}(A)$ such that $\mathbf{f}_n\rightarrow 0$ and $A\mathbf{f}_n\rightarrow \mathbf{h}$ with $\|\mathbf{h}\|_\mathcal{H}=1$. Since $A$ is dissipative, it follows that for every $\lambda>0$ and $\mathbf{f}\in \mathrm{Dom}(A)$
 \begin{eqnarray*}
\begin{aligned}
	\|(\frac{1}{\lambda}-A)(\mathbf{f}+\lambda^{-1}\mathbf{f}_n)\|\geq \frac{1}{\lambda}\|\mathbf{f}+\lambda^{-1}\mathbf{f}_n\|,
\end{aligned}
\end{eqnarray*}
that is,
% \begin{eqnarray*}
%\begin{aligned}
%	\|(\mathbf{f}+\lambda^{-1}\mathbf{f}_n)-\lambda A(\mathbf{f}+\lambda^{-1}\mathbf{f}_n)\|\geq\|\mathbf{f}+\lambda^{-1}\mathbf{f}_n\|,
%\end{aligned}
%\end{eqnarray*}
 \begin{eqnarray*}
\begin{aligned}
	\|(\mathbf{f}+\lambda^{-1}\mathbf{f}_n)-(\lambda A\mathbf{f}+ A\mathbf{f}_n)\|\geq\|\mathbf{f}+\lambda^{-1}\mathbf{f}_n\|.
\end{aligned}
\end{eqnarray*}
Letting $n\rightarrow\infty$ and then $\lambda\rightarrow 0$ gives $\|\mathbf{f}-\mathbf{h}\|\geq\|\mathbf{f}\|$, which is impossible obviously. Then we have $A$ is closable. 

From the inequality (\ref{dissipation}), we have for any $\mathbf{f}\in \mathrm{Dom}(A)$
\[
	\Re(A\mathbf{f}, \mathbf{f})_{\mathcal{H}}=-d\int_{\Gamma_1}|f_4|^2d\Gamma \leq 0,
\]
which implies that $A$ is a dissipative operator. For any $\mathbf{g}\in \mathrm{Dom}(A^{\ast})$, we have 
\begin{eqnarray*}
\begin{aligned}
	%\Re(A^{\ast}\mathbf{g}, \mathbf{g})_{\mathcal{H}^{\ast}}&=\Re(A^{\ast}\mathbf{g}, \mathbf{g})_{\mathcal{H}}\\
	&\Re(A^{\ast}\mathbf{g}, \mathbf{g})_{\mathcal{H}}\\
	&=\int_{D}\rho[c^{-2}(-c^2\Delta g_1)\overline{g_{2}} +\nabla g_2\overline{\nabla g_1}]dx+ \int_{\Gamma_{1}}[k(-g_3)\overline{g_{4}}+m(\frac{\rho}{m}g_{2} +\frac{k}{m}g_{3}  -\frac{d}{m}g_{4})\overline{g_{4}}]d\Gamma\\
	&=-d\int_{\Gamma_1}|g_4|^2d\Gamma \leq 0.
\end{aligned}
\end{eqnarray*}
Hence $A^{\ast}$ is also a dissipative operator.
Thus, from Corollary 4.4 in \cite[Chapter 1]{Pazy}, $A$ is the infinitesimal generator of a $C_0$ semigroup of contractions.

(ii) Thanks to $\frac{1}{2}\|X(t, \mathbf{x})\|^{2}_{\mathcal{H}}=\mathcal{E}(t)$, then from Young's inequality we have
 \begin{equation}\label{estimate1}
 \begin{aligned}
 	\|X(t, \mathbf{x})\|^{2}_{\mathcal{H}}&=\|\mathbf{x}\|^{2}_{\mathcal{H}}-2d\int_{0}^{t}\int_{\Gamma_1}\delta^{2}_{t}d\Gamma_1 dt-2\int_{\Gamma_1}\int_{0}^{t}\delta_{t}dW(s) d\Gamma_1 \\
	&\leq \|\mathbf{x}\|^{2}_{\mathcal{H}}-2d\int_{0}^{t}\int_{\Gamma_1}\delta^{2}_{t}d\Gamma_1 dt+2d\int_{\Gamma_1}\big(\int_{0}^{t}\delta_{t}dW(s)\big)^2 d\Gamma_1+\frac{|\Gamma_1|}{4d}.
\end{aligned}
\end{equation} 
By the It\^o isometry $\int_{\Gamma_1}\mathbb{E}\big(\int_{0}^{t}\delta_{t}dW(s)\big)^2 d\Gamma_1=\int_{\Gamma_1}\mathbb{E}\big(\int_{0}^{t}\delta^2_{t}ds\big) d\Gamma_1$,  it follows from (\ref{estimate1}) that
 \begin{equation}\label{estimate2}
 	\mathbb{E}\|X(t, \mathbf{x})\|^{2}_{\mathcal{H}}\leq\mathbb{E}\|\mathbf{x}\|^{2}_{\mathcal{H}}+\frac{|\Gamma_1|}{4d}.
\end{equation}
Since we have the relation $\sup\limits_{t>0}\mathbb{E}\|X(t, 0)\|^{2}_{\mathcal{H}}=\sup\limits_{t\geq0}\textrm{Tr}Q_t$, from (\ref{estimate2}) we get $\sup\limits_{t\geq0}\textrm{Tr}Q_t < \infty$. 

(iii) The existence and uniqueness of the mild solution to system (\ref{1}-\ref{2}-\ref{3}) is an immediate consequence of results (i) and (ii).
\end{proof}

Now, we give the statement and the proof of our main theorem as follows.
%It shows that the limit of the average of $f(X(t, \mathbf{x}))$ is a quantity that does not depend on the initial point. 

\begin{theorem} \label{mixing}
Suppose the boundary of the domain satisfies the geometric condition~\eqref{G1}. Then the system (\ref{1}-\ref{2}-\ref{3}) has a unique invariant measure $\mu$ satisfying
%\[
%	\lim_{t\rightarrow\infty}P(\mathbf{x}, \cdot) = \mu, \quad weakly
%\]
\[
	\lim_{t\rightarrow\infty}\mathbb{E}f(X(t, \mathbf{x}))=\int_{\mathcal{H}}f(z)\mu(dz)
\]
for any $f\in C_{b}(\mathcal{H})$, 
which implies the system (\ref{1}-\ref{2}-\ref{3}) is strongly mixing.
 
\end{theorem}
\begin{proof}
Based on Lemma \ref{operator}, there exists an invariant measure, denoted by $\mu$, for the system~(\ref{1}-\ref{2}-\ref{3}), which is derived from \cite[Theorem 6.2.1]{PratoZabczyk}. 
If the corresponding Markov semigroup $M_t$, $t\geqslant 0$, is $t_0$-regular for some $t_0>0$, then from the Doob's theorem~\cite{Doob} it follows that $\mu$ is the unique invariant measure for the semigroup $M_t$, which is also strongly mixing.  

It remains to show the regularity of the semigroup $M_t$. By Theorem \ref{regularity}, $M_t$ poccesses this regularity property if and only if there exists a constant $T_0>0$, depending only on $D$, such that for $T>T_0$, the following observability inequality holds
\begin{equation}\label{ob1}
	C(T)\|( \phi(0),  \phi_t(0), \delta(0), \delta_t(0))\|^2_{\mathcal{H}}\leq \int_{0}^{T}\int_{\Gamma_1}|\delta_t|^2 d\Gamma_1 dt
\end{equation}
for every solution of the adjoint system
\begin{equation} \label{ob}
\left
    \{
        \begin{array}{ll}
             \phi_{tt} = c^2\triangle \phi, \quad \textrm{in $[0, T]\times D$,}\\
            -\rho\phi_{t}= m\delta_{tt} - d\delta_t + k\delta, \quad  \frac{\partial\phi}{\partial \mathbf{n}}= \delta_{t}, \quad \textrm{on $[0, T]\times \Gamma_1$,}\\
            \phi = 0,  \quad \textrm{on $[0, T] \times \Gamma_0$,}   
        \end{array} 
\right.
\end{equation}
with any given value $( \phi(T),  \phi_t(T), \delta(T), \delta_t(T))$ at time $T$.
Next, we apply the multiplier method to prove the observability equality \eqref{ob1}. 

From Theorem A.4.1 in \cite{LTZ}, the geometric condition (\ref{G1}) implies that there exists a vector field $h(x)$ and a scalar function $H(x)$
\[
	h(x):=\nabla H(x) \in[C^2(\bar{D})]^3
\]
such that $\nabla H\cdot\mathbf{n}=0$ on $\Gamma_0$ and the Hessian matrix of $H$ evaluated on $\Gamma_0$ is positive definite, that is, $\nabla^2 H(x) \geqslant h_0 I$ for some constant $h_0 > 0$.
%Let $h$ be a $[C^{\infty}_{0}(\mathcal(R))]^{3}$-vector field, which will be given specifically later. 

Multiplying the first equation of (\ref{ob}) by $(h\cdot\nabla) \phi$, and integrating by parts, we have
\begin{equation}
\begin{aligned}\label{ob2}
	0&=\int_{\epsilon_0}^{T-\epsilon_0}\int_{D}(h\cdot\nabla) \phi ( \phi_{tt}  -c^2\triangle \phi)dxdt \\
	&=[\int_{D}(h\cdot\nabla) \phi  \phi_{t} dx]\bigg|_{\epsilon_0}^{T-\epsilon_0}-\int_{\epsilon_0}^{T-\epsilon_0}\int_{D}[\nabla\cdot(\frac{h}{2} \phi_{t}^2)-\frac{\nabla\cdot h}{2} \phi_{t}^2)]dxdt\\
	&\qquad-c^2\int_{\epsilon_0}^{T-\epsilon_0}\int_{\Gamma}(h\cdot\nabla) \phi\frac{\partial\phi}{\partial \mathbf{n}}d\Gamma dt+c^2\int_{\epsilon_0}^{T-\epsilon_0}\int_{D}\nabla(h\cdot\nabla \phi)\cdot \nabla\phi dxdt
\end{aligned}
\end{equation}
for some $\epsilon\in(0, \frac{T}{2})$. Since we have
\begin{eqnarray*}
\begin{aligned}
	\nabla(h\cdot\nabla \phi)\cdot \nabla\phi &= \sum_{i, j}\frac{\partial h_j}{\partial x_i}\frac{\partial\phi}{\partial x_i}\frac{\partial\phi}{\partial x_j} + \sum_{i, j}h_j\frac{\partial}{\partial x_j}(\frac{\partial\phi}{\partial x_i})^2\\
	&= \sum_{i, j}\frac{\partial h_j}{\partial x_i}\frac{\partial\phi}{\partial x_i}\frac{\partial\phi}{\partial x_j} + \nabla\cdot(\frac{h}{2}|\nabla\phi|^2)-\frac{\nabla\cdot h}{2}|\nabla\phi|^2,
\end{aligned}
\end{eqnarray*}
then from (\ref{ob2}) and using the boundary condition, it deduces that
\begin{eqnarray*}
\begin{aligned}
%	0&=[\int_{D}(h\cdot\nabla) \phi \frac{\partial \phi}{\partial t}dx]\bigg|_{\epsilon_0}^{T-\epsilon_0}-\int_{\epsilon_0}^{T-\epsilon_0}\int_{D}[\nabla\cdot(\frac{h}{2}(\frac{\partial \phi}{\partial t})^2)-c^2\nabla\cdot(\frac{h}{2}|\nabla\phi|^2)dxdt\\
%	&\qquad +\int_{\epsilon_0}^{T-\epsilon_0}\int_{D}[\frac{\nabla\cdot h}{2}(\frac{\partial \phi}{\partial t})^2)-c^2\frac{\nabla\cdot h}{2}|\nabla\phi|^2]dxdt\\
%	&\qquad-c^2\int_{\epsilon_0}^{T-\epsilon_0}\int_{\Gamma}(h\cdot\nabla) \phi\frac{\partial\phi}{\partial \mathbf{n}}d\Gamma dt+c^2\int_{\epsilon_0}^{T-\epsilon_0}\int_{D}\sum_{i, j}\frac{\partial h_j}{\partial x_i}\frac{\partial\phi}{\partial x_i}\frac{\partial\phi}{\partial x_j}dx dt
%	0&=\int_{\epsilon_0}^{T-\epsilon_0}\int_{D}\frac{\nabla\cdot h}{2}[(\frac{\partial \phi}{\partial t})^2-c^2|\nabla\phi|^2]dxdt-\int_{\epsilon_0}^{T-\epsilon_0}\int_{D}\nabla\cdot\{\frac{h}{2}[(\frac{\partial \phi}{\partial t})^2-c^2|\nabla\phi|^2]\}dxdt\\
%	&\qquad -c^2\int_{\epsilon_0}^{T-\epsilon_0}\int_{\Gamma}(h\cdot\nabla) \phi\frac{\partial\phi}{\partial \mathbf{n}}d\Gamma dt+c^2\int_{\epsilon_0}^{T-\epsilon_0}\int_{D}\sum_{i, j}\frac{\partial h_j}{\partial x_i}\frac{\partial\phi}{\partial x_i}\frac{\partial\phi}{\partial x_j}dx dt\\
%	&\qquad+[\int_{D}(h\cdot\nabla) \phi \frac{\partial \phi}{\partial t}dx]\bigg|_{\epsilon_0}^{T-\epsilon_0}\\
	0&=\int_{\epsilon_0}^{T-\epsilon_0}\int_{D}\frac{\nabla\cdot h}{2}[\phi_t^2-c^2|\nabla\phi|^2]dxdt-\int_{\epsilon_0}^{T-\epsilon_0}\int_{\Gamma}\{\frac{h\cdot\mathbf{n}}{2}[\phi_t^2-c^2|\nabla\phi|^2]\}d\Gamma dt\\
	&\qquad -c^2\int_{\epsilon_0}^{T-\epsilon_0}\int_{\Gamma_1}(h\cdot\nabla) \phi\delta_t d\Gamma_1 dt+c^2\int_{\epsilon_0}^{T-\epsilon_0}\int_{D}\sum_{i, j}\frac{\partial h_j}{\partial x_i}\frac{\partial\phi}{\partial x_i}\frac{\partial\phi}{\partial x_j}dx dt\\
	&\qquad+[\int_{D}(h\cdot\nabla) \phi \phi_t dx]\bigg|_{\epsilon_0}^{T-\epsilon_0}.
\end{aligned}
\end{eqnarray*}
Then we obtain
\begin{equation*}
\begin{aligned}
	&c^2 h_{0}\int_{\epsilon_0}^{T-\epsilon_0}\int_{D}|\nabla\phi|^2dxdt\\
	&\leq c^2\int_{\epsilon_0}^{T-\epsilon_0}\int_{D}\sum_{i, j}\frac{\partial h_j}{\partial x_i}\frac{\partial\phi}{\partial x_i}\frac{\partial\phi}{\partial x_j}dx dt\\
	&\leq \int_{\epsilon_0}^{T-\epsilon_0}\int_{\Gamma_1}\{\frac{h\cdot\mathbf{n}}{2}[\phi_t^2-c^2|\nabla\phi|^2]\}d\Gamma dt - \int_{\epsilon_0}^{T-\epsilon_0}\int_{D}\frac{\nabla\cdot h}{2}[\phi_t^2-c^2|\nabla\phi|^2]dxdt\\
	&\qquad +c^2\int_{\epsilon_0}^{T-\epsilon_0}\int_{\Gamma_1}(h\cdot\nabla) \phi\delta_t d\Gamma dt-[\int_{D}(h\cdot\nabla) \phi \phi_t dx]\bigg|_{\epsilon_0}^{T-\epsilon_0},
\end{aligned}
\end{equation*}
which implies
\begin{equation}
\begin{aligned} \label{ob3}
	c^2 h_{0}\int_{\epsilon_0}^{T-\epsilon_0}\int_{D}|\nabla\phi|^2dxdt\leq& C(h)\{\int_{\epsilon_0}^{T-\epsilon_0}\int_{\Gamma_1}[\phi_t^2 + |\nabla_{\|}\phi|^2]d\Gamma dt+\int_{0}^{T}\int_{\Gamma_1}\delta_t^2d\Gamma dt\}\\
	&+|\int_{\epsilon_0}^{T-\epsilon_0}\int_{D}\frac{\nabla\cdot h}{2}[\phi_t^2-c^2|\nabla\phi|^2]dxdt| +C(h)\mathcal{E}(T)
\end{aligned}
\end{equation}
due to $|\phi_t|=|\delta_t|^2 + |\nabla_{\|}\phi|^2$ on $\Gamma_1$ and $\mathcal{E}(T)=\mathcal{E}(s)+d\int_{s}^{T}\int_{\Gamma_1}\delta_{t}^2d\Gamma_1 dt$, $s\in[0, T]$.
Here, $\nabla_{\|}$ is the tangential derivative on the boundary and $\mathcal{E}(s) = \frac{1}{2}\|( \phi(s),  \phi_t(s), \delta(s), \delta_t(s))\|^2_{\mathcal{H}}$.

Multiplying the first equation of (\ref{ob}) by $\phi(\nabla\cdot h)$, and integrating by parts, we have
\begin{eqnarray*}
\begin{aligned}
	0&=\int_{\epsilon_0}^{T-\epsilon_0}\int_{D}\phi(\nabla\cdot h) (\phi_{tt}-c^2\triangle \phi)dxdt \\
	&=[\int_{D}\phi(\nabla\cdot h) \phi_{t}dx]\bigg|_{\epsilon_0}^{T-\epsilon_0}-\int_{\epsilon_0}^{T-\epsilon_0}\int_{D}[(\nabla\cdot h)\phi_t^2-c^2 \nabla( \phi \nabla\cdot h)\cdot \nabla\phi]dxdt\\
	&\qquad-c^2\int_{\epsilon_0}^{T-\epsilon_0}\int_{\Gamma} \phi(\nabla\cdot h)\frac{\partial\phi}{\partial \mathbf{n}}d\Gamma dt\\
	&=[\int_{D}\phi(\nabla\cdot h) \phi_{t}dx]\bigg|_{\epsilon_0}^{T-\epsilon_0}-\int_{\epsilon_0}^{T-\epsilon_0}\int_{D}(\nabla\cdot h)[\phi_t^2-c^2 |\nabla\phi|^2]dxdt\\
	&\qquad+c^2 \int_{\epsilon_0}^{T-\epsilon_0}\int_{D}  \phi \nabla\phi \cdot \nabla(\nabla\cdot h)dxdt-c^2\int_{\epsilon_0}^{T-\epsilon_0}\int_{\Gamma_1} \phi(\nabla\cdot h)\delta_t d\Gamma dt,
\end{aligned}
\end{eqnarray*}
which deduces that
\begin{equation}
\begin{aligned} \label{ob4}
	&\int_{\epsilon_0}^{T-\epsilon_0}\int_{D}(\nabla\cdot h)[\phi_t^2-c^2 |\nabla\phi|^2]dxdt\\
	&=\langle\phi(\nabla\cdot h), \phi_{t}\rangle\bigg|_{\epsilon_0}^{T-\epsilon_0}+c^2 \int_{\epsilon_0}^{T-\epsilon_0}\int_{D}  \phi \nabla\phi \cdot \nabla(\nabla\cdot h)dxdt-c^2\int_{\epsilon_0}^{T-\epsilon_0}\int_{\Gamma_1} \phi(\nabla\cdot h)\delta_t d\Gamma dt.	
\end{aligned}
\end{equation}
Here, $\langle\cdot, \cdot\rangle$ means the the product in $H^{-\eta}(D)\times H^{\eta}(D)$ for $\eta>0$. We introduce the standard denotation for the terms which are below the level of energy, that is,
\[
	\textrm{LOT}(\Phi, \Psi):=\|(\Phi, \Psi)\|^2_{C([0, T]; \mathcal{M})},
\]
in which $\Phi=\{\phi, \phi_t\}$, $\Psi=\{\delta, \delta_t\}$ and $\mathcal{M}=H^{1-\eta}(D)\times H^{-\eta}(D)\times H^{-\eta}(\Gamma_1)\times H^{-\eta}(\Gamma_1)$. Then by~(\ref{ob4}) and using the Young's inequality, we have
\begin{equation} \label{ob5}
\begin{aligned} 
	&\int_{\epsilon_0}^{T-\epsilon_0}\int_{D}(\nabla\cdot h)[\phi_t^2-c^2 |\nabla\phi|^2]dxdt\\
	& \leq \tau_1\int_{\epsilon_0}^{T-\epsilon_0}\int_{D}|\nabla\phi|^2dxdt+C(\tau_1) \int_{0}^{T}\int_{\Gamma_1}\delta_t^2d\Gamma dt +\textrm{LOT}(\Phi, \Psi).
\end{aligned}
\end{equation}
Combining (\ref{ob3}) and (\ref{ob5}), it deduces that
\begin{equation} 
\begin{aligned}\label{ob6}
	&\int_{\epsilon_0}^{T-\epsilon_0}\int_{D}\rho[c^{-2}|\phi_t|^2+ |\nabla\phi|^2]dxdt \\
	&\leq C(\tau_1, h) \{\int_{0}^{T}\int_{\Gamma_1}\delta_t^2d\Gamma dt+\int_{\epsilon_0}^{T-\epsilon_0}\int_{\Gamma_1}[\phi_t^2 + |\nabla_{\|}\phi|^2]d\Gamma dt\} +C(h)\mathcal{E}(T)+\textrm{LOT}(\Phi, \Psi).
\end{aligned}
\end{equation}

Let $\kappa(t)\in C_{0}^{\infty}(\mathbb{R})$ be cut-off function given by
\begin{equation*} 
\kappa(t)=\left
    \{
        \begin{array}{ll}
            1,\quad &t\in[\epsilon_0, T-\epsilon_0],\\
            \textrm{a $C^{\infty}$ function with range in $(0, 1)$}, \quad &t\in(0, \epsilon_0)\bigcup (T-\epsilon_0, T),\\
            0, \quad &t\in(-\infty, 0)\bigcup(T, \infty].    
        \end{array} 
\right.
\end{equation*}
%Estimate $\int_{\epsilon_0}^{T-\epsilon_0}\int_{\Gamma_1} |\nabla_{\|}\phi|^2d\Gamma_1dt$. 
Using Lemma 7.2 in \cite{LT}, we have
\begin{eqnarray*} 
\begin{aligned}
	&\int_{\epsilon_0}^{T-\epsilon_0}\int_{\Gamma_1} |\nabla_{\|}\phi|^2d\Gamma dt\\
	&\leq \int_{0}^{T}\int_{\Gamma_1} |\nabla_{\|}(\kappa\phi)|^2d\Gamma dt\\
	&\leq C(T, \epsilon_0)\{\int_{0}^{T}\int_{\Gamma_1} |\frac{\partial \phi}{\partial \mathbf{n}}|^2d\Gamma dt+\int_{0}^{T}\int_{\Gamma_1}| \frac{\partial}{\partial t}(\kappa \phi)|^2d\Gamma dt\}+\textrm{LOT}(\Phi, \Psi)\\
	&\leq C(T, \epsilon_0)\{\int_{0}^{T}\int_{\Gamma_1} |\delta_t|^2d\Gamma dt+\int_{0}^{T}\int_{\Gamma_1} \kappa^2|\phi_t|^2d\Gamma dt\}+\textrm{LOT}(\Phi, \Psi).
\end{aligned}
\end{eqnarray*}
We estimate $\int_{\epsilon_0}^{T-\epsilon_0}\int_{\Gamma_1}\phi_t^2 d\Gamma dt$ and $\int_{0}^{T}\int_{\Gamma_1} \kappa^2|\phi_t|^2d\Gamma dt$. The boundary condition on $\Gamma_1$ gives
\begin{eqnarray*} 
\begin{aligned}
	\int_{\epsilon_0}^{T-\epsilon_0}\int_{\Gamma_1}\phi_t^2 d\Gamma dt \leq \int_{0}^{T}\int_{\Gamma_1} \kappa^2|\phi_t|^2d\Gamma dt\leq C(\rho, m, d, k) \int_{0}^{T}\int_{\Gamma_1}[\delta^{2}_{tt}+\delta^{2}_{t}+\delta^{2}]d\Gamma dt.
\end{aligned}
\end{eqnarray*}
Then from (\ref{ob6}) it deduces that
%\begin{eqnarray*} 
%\begin{aligned}
%	\int_{\epsilon_0}^{T-\epsilon_0}\mathcal{E}(t)dt&=\frac{1}{2}\int_{\epsilon_0}^{T-\epsilon_0}\Big\{\int_{D}\rho[c^{-2}|\phi_t|^2+ |\nabla\phi|^2]dx+ \int_{\Gamma_{1}}[k|\delta|^{2}+m|\delta_t|^{2}]d\Gamma_1\Big\} dt\\
%	&\leq C(\tau, h, T, \epsilon_0, k, m)\int_{0}^{T}\int_{\Gamma_1}[\delta^{2}_{tt}+\delta^{2}_{t}+\delta^{2}]d\Gamma_1 dt+C(h)\mathcal{E}(0)+\textrm{LOT}(\Phi, \Psi).
%\end{aligned}
%\end{eqnarray*}
\begin{equation} 
\begin{aligned}\label{ob9}
	\int_{\epsilon_0}^{T-\epsilon_0}\mathcal{E}(t)dt \leq &C(\tau_1, h, T, \epsilon_0, \rho, m, d, k)\int_{0}^{T}\int_{\Gamma_1}[\delta^{2}_{tt}+\delta^{2}_{t}+\delta^{2}]d\Gamma dt\\
	&+C(h)\mathcal{E}(T)+\textrm{LOT}(\Phi, \Psi).
\end{aligned}
\end{equation}
From $\mathcal{E}(T)=\mathcal{E}(s)+d\int_{s}^{T}\int_{\Gamma_1}\delta_{t}^2d\Gamma_1 dt$, it follows that
\begin{eqnarray*} 
\begin{aligned}
	(T-2\epsilon_0)\mathcal{E}(0)&\leq (T-2\epsilon_0)\mathcal{E}(T)\\
	&= \int_{\epsilon_0}^{T-\epsilon_0}[ \mathcal{E}(s)+d\int_{s}^{T}\int_{\Gamma_1}\delta_{t}^2d\Gamma dt]ds\\	
	&\leq C(\tau_1, h, T, \epsilon_0, \rho, m, d, k)\int_{0}^{T}\int_{\Gamma_1}[\delta^{2}_{tt}+\delta^{2}_{t}+\delta^{2}]d\Gamma dt+C(h)\mathcal{E}(T)+\textrm{LOT}(\Phi, \Psi),
\end{aligned}
\end{eqnarray*}
which implies
\begin{equation} 
\begin{aligned}\label{ob10}
	(T-2\epsilon_0-C(h))\mathcal{E}(0)\leq C(\tau_1, h, T, \epsilon_0, \rho, m, d, k)\int_{0}^{T}\int_{\Gamma_1}[\delta^{2}_{tt}+\delta^{2}_{t}+\delta^{2}]d\Gamma dt+\textrm{LOT}(\Phi, \Psi).
\end{aligned}
\end{equation}

Differentiating the boundary equation in time, multiplying $\delta_t$ and integrating by parts, we have
\begin{eqnarray*} 
\begin{aligned}
	0&=\int_{0}^{T}\int_{\Gamma_1}\delta_{t}(\rho\phi_{tt}+m\delta_{ttt} - d\delta_{tt}+k\delta_{t})d\Gamma dt\\
	&=[\rho\int_{\Gamma_1}\delta_{t}\phi_{t}d\Gamma]\bigg|_{0}^{T}-\rho\int_{0}^{T}\int_{\Gamma_1}\delta_{tt}\phi_{t}d\Gamma dt+[m\int_{\Gamma_1}\delta_{t}\delta_{tt}d\Gamma_1]\bigg|_{0}^{T}-m\int_{0}^{T}\int_{\Gamma_1}\delta^{2}_{tt}d\Gamma dt\\
	&\qquad -[\frac{d}{2}\int_{\Gamma_1}\delta_{t}^2d\Gamma]\bigg|_{0}^{T}+k\int_{0}^{T}\int_{\Gamma_1}\delta^{2}_{t}d\Gamma dt,
\end{aligned}
\end{eqnarray*}
which implies
\begin{equation} 
\begin{aligned}\label{ob10-1}
	\int_{0}^{T}\int_{\Gamma_1}\delta^{2}_{tt}d\Gamma_1 dt%&=[\frac{\rho}{m}\int_{\Gamma_1}\delta_{t}\phi_{t}d\Gamma_1]\bigg|_{0}^{T}+[\int_{\Gamma_1}\delta_{t}\delta_{tt}d\Gamma_1]\bigg|_{0}^{T}+[\frac{d}{2m}\int_{\Gamma_1}\delta_{t}^2d\Gamma_1]\bigg|_{0}^{T}\\
	%&\qquad +\frac{k}{m}\int_{0}^{T}\int_{\Gamma_1}\delta^{2}_{t}d\Gamma_1 dt-\frac{\rho}{m}\int_{0}^{T}\int_{\Gamma_1}\delta_{tt}\phi_{t}d\Gamma_1 dt\\
	&\leq C(\rho, m, d)\mathcal{E}(0)+\frac{k}{m}\int_{0}^{T}\int_{\Gamma_1}\delta^{2}_{t}d\Gamma dt+\frac{\rho}{2m}\int_{0}^{T}\int_{\Gamma_1}(\delta_{tt}+\phi_{t})^2d\Gamma dt\\
	&\leq C(\rho, m, d)\mathcal{E}(0)+\frac{k}{m}\int_{0}^{T}\int_{\Gamma_1}\delta^{2}_{t}d\Gamma dt+\frac{\rho}{2m}\int_{0}^{T}\int_{\Gamma_1}(\delta_{t}+\delta)^2d\Gamma dt\\
	&\leq C(\rho, m, d)\mathcal{E}(0)+\frac{k+\rho}{m}\int_{0}^{T}\int_{\Gamma_1}\delta^{2}_{t}d\Gamma dt+\frac{\rho}{m}\int_{0}^{T}\int_{\Gamma_1}\delta^2d\Gamma dt.
\end{aligned}
\end{equation}
Multiplying $\delta$ and integrating by parts, we have
\begin{eqnarray*} 
\begin{aligned}
	0&=\int_{0}^{T}\int_{\Gamma_1}\delta(\rho\phi_{t}+m\delta_{tt}-d\delta_{t}+k\delta)d\Gamma dt\\
	&=[\rho\int_{\Gamma_1}\delta\phi d\Gamma]\bigg|_{0}^{T}-\rho\int_{0}^{T}\int_{\Gamma_1}\delta_{t}\phi d\Gamma dt+[m\int_{\Gamma_1}\delta\delta_{t}d\Gamma_1]\bigg|_{0}^{T}-m\int_{0}^{T}\int_{\Gamma_1}\delta^{2}_{t}d\Gamma dt\\
	&\qquad -[\frac{d}{2}\int_{\Gamma_1}\delta^2d\Gamma]\bigg|_{0}^{T}+k\int_{0}^{T}\int_{\Gamma_1}\delta^{2}d\Gamma dt,
\end{aligned}
\end{eqnarray*}
which implies
\begin{equation} \label{ob10-2}
\begin{aligned}
	\int_{0}^{T}\int_{\Gamma_1}\delta^{2}d\Gamma dt%&=[-\frac{\rho}{k}\int_{\Gamma_1}\delta\phi d\Gamma_1]\bigg|_{0}^{T}-[\frac{m}{k}\int_{\Gamma_1}\delta\delta_{t}d\Gamma_1]\bigg|_{0}^{T}-[\frac{d}{2k}\int_{\Gamma_1}\delta^2d\Gamma_1]\bigg|_{0}^{T}\\
	%&\qquad +\frac{\rho}{k}\int_{0}^{T}\int_{\Gamma_1}\delta_{t}\phi d\Gamma_1 dt+\frac{m}{k}\int_{0}^{T}\int_{\Gamma_1}\delta^{2}_{t}d\Gamma_1 dt\\
	&\leq C(\rho, m, d, k)\mathcal{E}(0)+\tau_2\int_{0}^{T}\int_{\Gamma_1}\phi^2 d\Gamma dt+(C(\tau_2, \rho)+m)\int_{0}^{T}\int_{\Gamma_1}\delta_{t}^2 d\Gamma dt\\
	&\leq [C(\rho, m, d, k)+\tau_2 T]\mathcal{E}(0)+(C(\tau_2, \rho)+\frac{m}{k})\int_{0}^{T}\int_{\Gamma_1}\delta_{t}^2 d\Gamma dt. 
\end{aligned}
\end{equation}
Then from (\ref{ob10}) and choosing $T_{0}:=\frac{2\epsilon_0 + C(h) + C(\rho, m, d, k)}{1-\tau_2}$, combining (\ref{ob10-1}) and (\ref{ob10-2}) gives
%\begin{eqnarray} 
%\begin{aligned}\label{ob13}
%	&[T-2\epsilon_0-C(h)-C(\rho, m, d, k)-\tau_2 T]\mathcal{E}(0)\\
%	&\leq C(\tau, \tau_2, h, T, \epsilon_0, \rho, m, d, k)\int_{0}^{T}\int_{\Gamma_1}\delta^{2}_{t}d\Gamma_1 dt+\textrm{LOT}(\Phi, \Psi),
%\end{aligned}
%\end{eqnarray}
\begin{equation} 
\begin{aligned}\label{ob13}
	(T-T_0)\mathcal{E}(0)
	\leq C(\tau_1, \tau_2, h, T, \epsilon_0, \rho, m, d, k)\int_{0}^{T}\int_{\Gamma_1}\delta^{2}_{t}d\Gamma dt+\textrm{LOT}(\Phi, \Psi).
\end{aligned}
\end{equation}

Claim that there exists a constant $C_T>0$ such that the solution of (\ref{ob}) satisfies the inequality
\begin{equation}\label{ob14}
	\textrm{LOT}(\Phi, \Psi)\leq C_T \|\delta_{t}\|^2_{L^2([0, T]\times \Gamma_1)}.
\end{equation}
We use the method of contradiction. Suppose that this claim is false. Then there exists a sequence of solutions $\{\Phi^{(n)}(t), \Psi^{(n)}(t)\}$ to (\ref{ob}) with the initial data 
$\{\Phi^{(n)}(0), \Psi^{(n)}(0)\}\subset\mathcal{H}$ such that
\begin{equation} \label{false}
\begin{aligned}
	\textrm{LOT}(\Phi^{(n)}, \Psi^{(n)}) =1, \quad \|\delta^{(n)}_{t}\|^2_{L^2([0, T]\times \Gamma_1)}\xrightarrow{n\rightarrow \infty} 0.
\end{aligned}
\end{equation}
From (\ref{ob13}), it follows that $\|\{\Phi^{(n)}(0), \Psi^{(n)}(0)\}\|_{\mathcal{H}}$ is bounded. Then there exists a subsequence, still denoted by $\{\Phi^{(n)}(0), \Psi^{(n)}(0)\}$, and $(\hat{\Phi}(0), \hat{\Psi}(0))\in\mathcal{H}$ such that
\begin{eqnarray*}
	\{\Phi^{(n)}(0), \Psi^{(n)}(0)\}\xrightarrow{weakly} \{\hat{\Phi}(0), \hat{\Psi}(0)\},\quad \textrm{in $\mathcal{H}$}
\end{eqnarray*}
as $n\rightarrow\infty$. Let $\{\hat{\Phi}(t), \hat{\Psi}(t)\}$ be the solution of (\ref{ob}) subject to the initial data $\{\hat{\Phi}(0), \hat{\Psi}(0)\}$. By~(\ref{dissipation}), we see that $\|\{\Phi^{(n)}(t), \Psi^{(n)}(t)\}\|_{C([0, T]; \mathcal{H})}$ is bounded, and then
\begin{equation}\label{ob15}
	\{\Phi^{(n)}(t), \Psi^{(n)}(t)\}\xrightarrow{\textrm{weak star}} \{\hat{\Phi}(t), \hat{\Psi}(t)\},\quad \textrm{in $L^{\infty}([0, T]; \mathcal{H})$}.
\end{equation}
Let $\mathcal{N}=H^{-\eta}(D)\times (H^{1}(D))^{\prime}\times H^{-\eta}(\Gamma_1)\times H^{-\eta}(\Gamma_1)$. For any $\omega\in H^1$ and $t\in(0, T)$, we have the following equation holds true
\[
 (\phi^{(n)}_{tt}, \omega)_{H^{-1}\times H^1} =\int_{D}\Delta \phi^{(n)} \omega dx = \int_{\Gamma_1}\delta^{(n)}_t\omega d\Gamma_1-\int_{D}\nabla \phi^{(n)}\cdot \nabla\omega dx, 
 \]
 which implies that  
 \[
 	|(\phi^{(n)}_{tt}, \omega)_{H^{-1}\times H^1}|\leq C(\|\delta^{(n)}\|_{L^2}, \|\nabla\phi^{(n)}\|_{L^2})\|\omega\|_{H^1}.
\]
Hence, $\phi^{(n)}_{tt}\in L^{\infty}([0, T], (H^{1}(D))^{\prime})$. Thus, it follows that $\|\{\Phi^{(n)}(t), \Psi^{(n)}(t)\}\|_{L^{\infty}([0, T], \mathcal{N})}$ is bounded uniformly. Due to $\mathcal{H}\subset\subset\mathcal{M}\subset\mathcal{N}$, we deduce from the Aubin's theorem \cite{Simon} that
 \begin{eqnarray*}
	\{\Phi^{(n)}(t), \Psi^{(n)}(t)\}\xrightarrow{\textrm{strongly}} \{\hat{\Phi}(t), \hat{\Psi}(t)\},\quad \textrm{in $L^{\infty}([0, T]; \mathcal{M})$},
\end{eqnarray*} 
which gives
\begin{equation}\label{ob16}
\|\{\hat{\Phi}(t), \hat{\Psi}(t)\}\|_{C([0, T]; \mathcal{M})}=\textrm{LOT}(\Phi^{(n)}, \Psi^{(n)}) =1.
\end{equation} 
By (\ref{false}) and (\ref{ob15}), we have $\hat{\delta}_t =0$. Let $P=\hat{\phi}_{tt}$. Then we have
\begin{equation*}
\left
    \{
        \begin{array}{ll}
             P_{tt} = c^2\triangle P, \quad \textrm{in $\mathbb{R}^{+} \times D$,}\\
            \frac{\partial P}{\partial \mathbf{n}}= 0, \quad \textrm{on $\mathbb{R}^{+} \times \Gamma_1$,}\\
            P = 0,  \quad \textrm{on $\mathbb{R}^{+} \times \Gamma_0$}.   
        \end{array} 
\right.
\end{equation*}
By Holmgren's Uniqueness Theorem, choosing $T>2\textrm{diam}(D)$, it follows that $P=\hat{\phi}_{tt}=0$. And then we obtain
\begin{equation*}
\left
    \{
        \begin{array}{ll}
             \triangle\hat{\phi} = 0, \quad \textrm{in $ D$,}\\
            \frac{\partial \hat{\phi}}{\partial \mathbf{n}}= 0, \quad \textrm{on $ \Gamma_1$,}\\
            \hat{\phi}= 0,  \quad \textrm{on $ \Gamma_0$},  
        \end{array} 
\right.
\end{equation*}
which implies that $\hat{\phi}=0$. The boundary condition shows $\hat{\delta}=0$. Then we obtain $\{\hat{\Phi}(t), \hat{\Psi}(t)\}=0$ which contradicts (\ref{ob16}). Therefore, we prove the inequality (\ref{ob14}).  

In summary, the inequality \eqref{ob1} follows from \eqref{ob13} and \eqref{ob14}.
\end{proof}

\section{Conclusion}
\label{sec:conclusion}
This paper is the first part of the project devoted to studying the problem of mixing for acoustic wave motion driven by a boundary random force. 
We need to point out that the acoustic boundary condition imposed on the portion $\Gamma_1$ of the boundary is locally reacting in the present paper.  
However, the motion along $\Gamma_1$ can be modeled as a vibrating membrane, which propagates waves along its surface with the speed $c^2_{\Gamma_1}$. Such boundary conditions obtained below are called non-locally reacting. 
\begin{equation*} %\tag{B2-2}
\left
    \{
        \begin{array}{ll}
            -\rho\phi_{t}= m\delta_{tt} -c^2_{\Gamma_1}\triangle_{\Gamma_1}\delta+d\delta_t+k\delta + \xi(t),\\ 
            \frac{\partial\phi}{\partial \mathbf{n}}= \delta_{t},  \label{4}
        \end{array} 
\right. \quad \textrm{on $\mathbb{R}^{+} \times \Gamma_1$}.
\end{equation*}
As for $\xi = 0$, a good literature review can be found in \cite{mugnolo2021}. For the case $\xi\neq 0$, the investigation of the mixing problem for acoustic wave system with the above non-locally reacting boundary condition as the second part of our project is ongoing.
Moreover, it is interesting and applicable to consider the mixing problem for acoustic wave system with a nonlinear boundary condition which is our future research work.

\section*{Declaration of competing interest}
The authors declare that they have no known competing financial interests or personal relationships that could have appeared to influence the work reported in this paper.

\section*{References}

\bibliography{mybibfile}

\end{document}